\documentclass[11pt]{amsart}
\usepackage{latexsym,amssymb,amsmath}
\usepackage[all,ps,cmtip]{xy} 
\newtheorem{theo}{Theorem}

\newtheorem{prop}{Proposition}
\newtheorem{lemm}{Lemma}

\begin{document}

\def\ot{\otimes}
\def\we{\wedge}
\def\wec{\wedge\cdots\wedge}
\def\op{\oplus}
\def\ra{\rightarrow}
\def\lra{\longrightarrow}
\def\fso{\mathfrak so}
\def\cO{\mathcal{O}}
\def\cE{\mathcal{E}}\def\cF{\mathcal{F}}
\def\cS{\mathcal{S}}
\def\cL{\mathcal{L}}
\def\fsl{\mathfrak sl}\def\fe{\mathfrak e}
\def\PP{\mathbb P}\def\PP{\mathbb P}\def\ZZ{\mathbb Z}\def\CC{\mathbb C}
\def\RR{\mathbb R}\def\HH{\mathbb H}\def\OO{\mathbb O}\def\OP2{\OO\PP^2}
\def\smc{\cdots }\def\JO{{\mathcal J}_3(\OO)}
\def\tr{{\rm trace}\;}\def\a{\alpha}\def\om{\omega}

\title{On the derived category of the Cayley plane}

\author{L. Manivel}

\begin{abstract} 
We describe a maximal exceptional collection on the Cayley plane, the minimal homogeneous 
projective variety of $E_6$. This collection consists
in a sequence of $27$ irreducible homogeneous bundles. 
\end{abstract}

\maketitle
\section{The Cayley plane}

Let ${\bf O}$ denote the normed algebra of (real) octonions (see e.g. \cite{baez}), and let 
$\OO$ be its complexification. The space
$$\JO = \Bigg\{ 
\begin{pmatrix} 
c_1            & x_3            & \bar{x}_2  \\ 
\bar{x}_3 & c_2            &  x_1            \\ 
x_2            & \bar{x}_1 &  c_3            \\ 
\end{pmatrix}, 
\;\; c_i \in {\bf C}, \; x_i \in \OO \Bigg\} \cong {\bf C}^{27}
$$
of $\OO$-Hermitian matrices of order $3$, is the exceptional 
simple complex Jordan algebra.

The subgroup $SL_3(\OO)$ of $GL(\JO)$ consisting of automorphisms preserving the 
determinant is the adjoint group of type $E_6$. 
The action of $E_6$ on the 
projectivization $\PP\JO$ has exactly three orbits: the complement of the 
determinantal hypersurface, the regular part of this hypersurface, and its singular
part which is the closed $E_6$-orbit. These three orbits can be viewed 
as the (projectivized) sets of matrices of rank three, two, and one respectively. 

The closed orbit, i.e. the (projectivization of) the set of rank one matrices,
is called the {\it Cayley plane} and denoted $\OP2$. It can be defined by the quadratic equation 
$$X^2=\tr (X)X, \qquad  X\in \JO,$$ 
or as the closure of the affine cell 
$$\OO\PP^2_0=
\Bigg\{\begin{pmatrix} 1 & x & y \\ \bar{x} & x\bar{x} & y\bar{x} \\ 
\bar{y} & x\bar{y} & y\bar{y} \end{pmatrix}, \quad x,y\in\OO\Bigg\}\cong {\bf C}^{16}.$$

\medskip
Since the Cayley plane is a closed orbit of $E_6$, it can also be identified with the 
quotient of $E_6$ by a parabolic subgroup, namely the maximal parabolic subgroup $P_1$ defined 
by the simple root $\a_1$ in the notation below. The semi-simple part 
of this maximal parabolic is isomorphic to ${\rm Spin}_{10}$.
\medskip

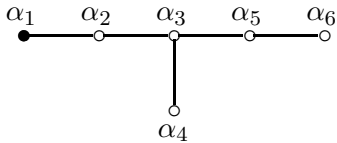
\begin{center}
\setlength{\unitlength}{5mm}
\begin{picture}(20,5)(-6.5,0)
\multiput(-.3,3.8)(2,0){5}{$\circ$}
\multiput(0,4)(2,0){4}{\line(1,0){1.7}}
\put(3.7,1.8){$\circ$}
\put(-.3,3.8){$\bullet$}
\put(-.6,4.4){$\a_1$}
\put(1.4,4.4){$\a_2$}
\put(3.4,4.4){$\a_3$}
\put(5.4,4.4){$\a_5$}
\put(7.4,4.4){$\a_6$}
\put(3.45,1.3){$\a_4$}
\put(3.89,2.15){\line(0,1){1.73}}
\end{picture}
\end{center}

If we denote by $V_\om$ the irreducible $E_6$-module with highest weight $\om$, we have $\JO\simeq V_{\om_1}$.
This  is a {\it minuscule} module, meaning that its weights with respect 
to any maximal torus of $E_6$, are all conjugate under the action of the Weyl group $W(E_6)$. 
For more details, see \cite{LM,IM}. 

Note that the Dynkin diagram of type $E_6$ has an obvious symmetry of order two, which accounts for 
the duality between irreducible modules. For example, the dual module of $V_{\om_1}$ is  $V_{\om_6}$.

\section{Homogeneous bundles on the Cayley plane}

\subsection{Irreducible homogeneous bundles}
The category of homogeneous bundles on a rational homogeneous variety $G/P$
is equivalent to the category $Mod_P$ of $P$-modules. Recall that $P$ has a 
non trivial decomposition $P=LP^u$, where $P^u$ denotes the unipotent radical
and $L$ a Levi factor. Since $P^u$ is non trivial, $P$ is not reductive,
and $P$-modules are not completely reducible in general. Irreducible $P$-modules 
have a trivial action of $P^u$, so that they are completely determined by 
their $L$-module structure. Since $L$ is reductive, its irreducible modules 
are well understood: they are uniquely determined by their highest weight $\om$, 
which can be any $L$-dominant weight of $G$. We denote by $\cE_\om$ the corresponding 
irreducible homogeneous vector bundles on $G/P$. By the Borel-Weil theorem, 
$H^0(G/P,\cE_\om)=V_\om^\vee$
if $\om$ is dominant, and otherwise $H^0(G/P,\cE_\om)=0$. 

\smallskip
For the Cayley plane $\OP2 = E_6/P_1$, a Levi factor $L$ of $P_1$, modded out by its
one dimensional center,  is a copy  
of $\mathrm{Spin}_{10}$. An $L$-dominant weight $\omega$ is a linear combination
$\om=a_1\om_1+a_2\om_2+a_3\om_3+a_4\om_4+a_5\om_5+a_6\om_6$ of the fundamental 
weights of $\fe_6$, with $a_2,\ldots ,a_6\ge 0$.  We can encode $\om$ by the
Dynkin diagram of $E_6$, where the node corresponding to the fundamental weight 
$\omega_i$ is labeled $a_i$. 

\medskip\noindent {\it Example 1}. The weight $\omega=-\omega_1$ defines 
a character of $L$. So  $\cE_{-\om_1}$ is just a line bundle, the negative 
generator of the Picard group.  The dual bundle  $\cE_{\om_1}$ defines the 
embedding of $\OP2$ in $\PP V_{\om_1}=\PP\JO$ and will be denoted $\cO_{\OP2}(1)$.  

\medskip
\begin{center}
\setlength{\unitlength}{5mm}
\begin{picture}(20,4)(-6.5,1)
\put(-6,2.8){$\cO_{\OP2}(1)\quad\simeq$}
\multiput(-.3,3.8)(2,0){5}{$\circ$}
\multiput(0,4)(2,0){4}{\line(1,0){1.7}}
\put(3.7,1.8){$\circ$}
\put(-.3,3.8){$\bullet$}
\put(-.3,4.5){$1$}
\put(3.89,2.15){\line(0,1){1.73}}
\end{picture}
\end{center}

\noindent {\it Example 2}. The weight that defines the tangent bundle of $\OP2$ is 
the highest root of $\fe_6$, which is also the dominant weight defining the 
adjoint representation. Note that the corresponding representation of $\mathrm{Spin}_{10}$
is one of the half-spin representations, which has dimension sixteen, as the Cayley plane. 

\smallskip
\begin{center}
\setlength{\unitlength}{5mm}
\begin{picture}(20,4)(-6.5,.5)
\put(-5,2.8){$T_{\OP2}\quad\simeq$}
\multiput(-.3,3.8)(2,0){5}{$\circ$}
\multiput(0,4)(2,0){4}{\line(1,0){1.7}}
\put(3.7,1.8){$\circ$}
\put(-.3,3.8){$\bullet$}
\put(3.7,1){$1$}
\put(3.89,2.15){\line(0,1){1.73}}
\end{picture}
\end{center}

The Borel-Weil theorem yields that $H^0(\OP2,T_{\OP2})=\fe_6$, as expected. 

Since the two half-spin representations of $\mathrm{Spin}_{10}$ are dual one of the other,
one could expect that the weight defining the cotangent bundle of $\OP2$ be
$\om_2$. This is not exactly true: the defining weight is $\om_2-\om_1$, where 
substracting $\omega_1$ amounts, at the level of bundles, to twisting by
$\cO_{\OP2}(-1)$. To check this, one needs to remember that 
if an irreducible $L$-module has highest weight $\om$, then its lowest 
weight is $w_0^L(\om)$, where $w_0^L$ denotes the longest element of the 
Weyl group $W(L)$ of $L\simeq \mathrm{Spin}_{10}\times\CC^*$, and then the highest weight of the 
dual module is $-w_0^L(\om)$. But this weight must be computed 
inside  the weight lattice of $\fe_6$, on which $W(L)$ acts naturally since it is
embedded in $W(E_6)$. And the result of this computation will be what it would 
be in the weight lattice of $\mathrm{Spin}_{10}$, only up to extra multiples of $\om_1$.

\medskip
\begin{center}
\setlength{\unitlength}{5mm}
\begin{picture}(20,5)(-6.5,0)
\put(-5,2.8){$\Omega^1_{\OP2}\quad\simeq$}
\multiput(-.3,3.8)(2,0){5}{$\circ$}
\multiput(0,4)(2,0){4}{\line(1,0){1.7}}
\put(3.7,1.8){$\circ$}
\put(-.3,3.8){$\bullet$}
\put(-.6,4.4){$-1$}
\put(1.7,4.4){$1$}
\put(3.89,2.15){\line(0,1){1.73}}
\end{picture}
\end{center}

\medskip\noindent {\it Example 3}. The minimal non trivial representation 
of $\mathrm{Spin}_{10}$ is the vector representation. This implies
that up to line bundles, the irreducible homogeneous 
bundle defined by $\om_6$ has minimal rank, equal to ten. 
We denote it by $\cS$.

\medskip
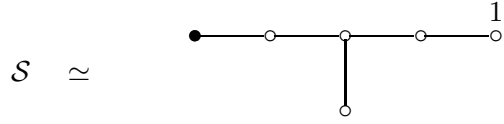
\begin{center}
\setlength{\unitlength}{5mm}
\begin{picture}(20,4)(-6.5,1)
\put(-5,2.8){$\cS\quad\simeq$}
\multiput(-.3,3.8)(2,0){5}{$\circ$}
\multiput(0,4)(2,0){4}{\line(1,0){1.7}}
\put(3.7,1.8){$\circ$}
\put(-.3,3.8){$\bullet$}
\put(7.7,4.4){$1$}
\put(3.89,2.15){\line(0,1){1.73}}
\end{picture}
\end{center}

The vector representation of $\mathrm{Spin}_{10}$ is self-dual. For the reasons explained above, this does 
not quite imply that $\cS$ be self-dual, but this must be the case up to a twist by a line bundle. 
One can easily check that $\cS^\vee=\cS(-1)$.
 
The geometric interpretation of $\cS$ is the following. By the Borel-Weil
theorem, we have  $H^0(\OO\PP^2,\cS)=V_{\omega_6}^\vee=V_{\omega_1}=\JO$. An irreducible 
homogeneous bundle with non trivial sections is generated by global sections, 
so dualizing the evaluation map we get an injection 
$$\cS^\vee\hookrightarrow \JO^\vee\otimes\cO_{\OO\PP^2}.$$
This map identifies each fiber of $\cS^\vee$ with the linear span of an $\OO$-line,
a maximal quadric in the dual Cayley plane $\OO\PP^2\subset \PP \JO^\vee$, see \cite{LM}. 
(Note that 
the Cayley plane and its dual are isomorphic, but only non-canonically: this 
reflects the fact that the order two symmetry of  the Dynkin diagram can only be realized as 
an outer automorphism of $E_6$.) 
In particular the presence of this maximal quadric explains that there is a natural 
quadratic form $$Sym^2\cS\rightarrow \cO_{\OO\PP^2}(1).$$ 

\medskip\noindent  {\it Definition}.  
Let $\cS_2$ be the kernel of the map $Sym^2\cS\rightarrow \cO_{\OO\PP^2}(1)$. Since
the symmetric square of the vector representation of $\mathrm{Spin}_{10}$ is, up to the trivial
factor defined by the quadratic form, irreducible, $\cS_2$ is an irreducible vector bundle, 
with highest weight $2\om_6$. 

The quadratic map $Sym^2\cS\rightarrow \cO_{\OO\PP^2}(1)$  induces a cubic map   
$Sym^3 \cS\rightarrow \cS(1)$. Let $\cS_3$ be the kernel of this cubic map. 
This is the irreducible bundle with highest weight $3\om_6$.

\subsection{Bott's theorem}
The fundamental tool for computing the cohomology of vector bundles on homogeneous spaces 
is Bott's theorem, which extends the Borel-Weil theorem for global sections to higher cohomology 
groups. 

Consider on $G/P$ an irreducible vector bundle $\cE_\om$. We have seen that it has 
non trivial global sections exactly when $\om$ is dominant. In general, let $\rho$ 
denote the sum of the fundamental weights, and consider the weight $\om+\rho$. 
This weight is {\it singular} if there exists a root $\alpha$ such that $\langle 
\om+\rho,\alpha^\vee\rangle =0$ (equivalently, $\om+\rho$ is fixed by the simple reflection 
$\alpha$). Otherwise, there exists a unique $w$ in the Weyl group such that $w(\om+\rho)$
be strictly dominant, and then $w(\om+\rho)-\rho$ is dominant. 

\begin{theo}[Bott's theorem]
If $\om+\rho$ is singular, then $\cE_\om$ is acyclic. Otherwise, there is a unique $w\in W(E_6)$ such that 
$w(\om+\rho)$ is strictly dominant. Then 
$$H^{\ell(w)}(G/B,\cE_\om)=V_{w(\om+\rho)-\rho}^\vee,$$
and the other cohomology groups of $\cE_\om$ vanish.
\end{theo}

\noindent {\it Remark}. 
To check whether the weight $\omega+\rho$ is singular or not, we can proceed as follows. 
If $\omega+\rho$ is not dominant, one of its coefficients on the basis of fundamental weights, say on
$\om_i$,  must be negative. 
Then we apply the simple reflection $s_{\alpha_i}$, in order to cross the hyperplane orthogonal
to $\alpha_i^\vee$. Not that since $E_6$ is simply laced, this simply amounts to changing the 
(negative) coefficient of $\om_i$ into its opposite, and adding it to the coefficients of the 
fundamental weights connected to $\om_i$ in the Dynkin diagram. Iterating this procedure, 
we will eventually get a weight with some zero coefficient, which will imply  that $\om+\rho$ 
is singular, or get a strictly dominant weight which will be the representative $w(\om+\rho)$ 
of the $W(E_6)$-orbit of $\om+\rho$ in the interior of  the dominant chamber. In the latter case, 
the number of applications of these simple reflections is nothing but the length $\ell(w)$ of $w$, 
which is the degree of the only non-zero cohomology group of $\cE_\om$. 

\section{Exceptional sequences}

\subsection{Exceptional bundles}

Recall that an object $\cF$ of the derived category of coherent sheaves on a variety 
$X$ is {\it exceptional} if $RHom(\cF,\cF)=\CC$. If $\cF$ is represented by a vector bundle
$F$ on $X$, this means that 
$$H^i(X, End(F))=\delta_{i,0}\CC.$$
   
\begin{prop}
The homogeneous bundles $\cS, \cS_2, \cS_3$ on $\OO\PP^2$ are exceptional. 
\end{prop}

\proof If $U$ denotes the vector representation of $\mathrm{Spin}_{10}$, we know that 
$\wedge^2U$ is an irreducible (and even fundamental) representation, and that 
$Sym^2U$ splits into a trivial factor generated by the invariant quadratic form, 
and an irreducible summand. At the level of bundles, since $\cS^\vee=\cS(-1)$, 
this implies that 
$$End(\cS)=\cE_{\omega_5}(-1)\oplus \cO_{\OP2}\oplus \cS_2(-1).$$
The bundle $\cE_{\omega_5}(-1)$ has highest weight $\om=\om_5-\om_1$. Since
$\om+\rho=\om_2+\om_3+\om_4+2\om_5+\om_6$ is orthogonal to $\alpha_1^\vee$,
$\om+\rho$ is singular. By Bott's theorem we conclude that $\cE_{\omega_5}(-1)$
is acyclic. For exactly the same reason $\cS_2(-1)$ is also acyclic. We 
conclude that 
$$H^i(\OP2, End(\cS))=H^i(\OP2,\cO_{\OP2})=\delta_{i,0}\CC.$$
So $\cS$ is exceptional. 

We proceed similarly with the other two bundles. First observe that 
$\cS_2^\vee=\cS_2(-2)$ and $\cS_3^\vee=\cS_3(-3)$. Using e.g. LiE \cite{lie}
to compute tensor products of representations of $\mathrm{Spin}_{10}$, we get 
the decompositions:
\begin{eqnarray*}
End(\cS_2) &= & \cE_{4\omega_6}(-2)\oplus \cE_{\omega_5+2\omega_6}(-2)\oplus \cE_{2\omega_5}(-2)
\oplus \cE_{2\omega_6}(-1)\oplus \cE_{\omega_5}(-1)\oplus \cO_{\OP2}, \\
End(\cS_3) &= & \cE_{6\omega_6}(-3)\oplus \cE_{\omega_5+4\omega_6}(-3)\oplus \cE_{2\omega_5+2\omega_6}(-3)
\oplus \cE_{3\omega_5}(-3)\oplus \cE_{4\omega_6}(-2)\oplus  \\
 & &\oplus\; \cE_{\omega_5+2\omega_6}(-2)\oplus 
\cE_{\omega_3+2\omega_5}(-3)\oplus \cE_{2\omega_6}(-1)\oplus \cE_{\omega_5}(-1)\oplus \cO_{\OP2}.
\end{eqnarray*}
Our claim amounts to the acyclicity of all the non trivial vector bundles in these decompositions, 
hence, by Bott's theorem, to the singularity of all the corresponding weights, once we have added 
$\rho$. We use the remark after Bott's theorem above. Consider for example $\cE_{6\omega_6}(-3)$,
whose highest weight is $6\omega_6-3\omega_1$. After adding $\rho$,  we get successively, applying 
$s_{\alpha_1}$ and $s_{\alpha_2}$:

\medskip
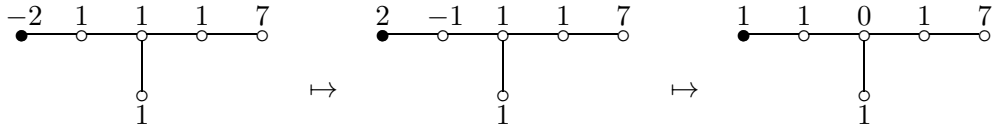
\begin{center}
\setlength{\unitlength}{4mm}
\begin{picture}(20,5)(5,0)
\multiput(-.3,3.8)(2,0){5}{$\circ$}
\multiput(0.1,4.08)(2,0){4}{\line(1,0){1.65}}
\put(3.7,1.8){$\circ$}
\put(-.3,3.8){$\bullet$}
\put(-.6,4.4){$-2$}
\put(1.7,4.4){$1$}
\put(3.7,4.4){$1$}
\put(5.7,4.4){$1$}
\put(7.7,4.4){$7$}
\put(3.7,1.1){$1$}
\put(3.94,2.18){\line(0,1){1.68}}

\put(9.5,2){$\mapsto$}

\multiput(11.7,3.8)(2,0){5}{$\circ$}
\multiput(12.1,4.08)(2,0){4}{\line(1,0){1.65}}
\put(15.7,1.8){$\circ$}
\put(11.7,3.8){$\bullet$}
\put(11.7,4.4){$2$}
\put(13.4,4.4){$-1$}
\put(15.7,4.4){$1$}
\put(17.7,4.4){$1$}
\put(19.7,4.4){$7$}
\put(15.7,1.1){$1$}
\put(15.94,2.18){\line(0,1){1.68}}

\put(21.5,2){$\mapsto$}

\multiput(23.7,3.8)(2,0){5}{$\circ$}
\multiput(24.1,4.08)(2,0){4}{\line(1,0){1.65}}
\put(27.7,1.8){$\circ$}
\put(23.7,3.8){$\bullet$}
\put(23.7,4.4){$1$}
\put(25.7,4.4){$1$}
\put(27.7,4.4){$0$}
\put(29.7,4.4){$1$}
\put(31.7,4.4){$7$}
\put(27.7,1.1){$1$}
\put(27.94,2.18){\line(0,1){1.68}}
\end{picture}
\end{center}

Since there is a zero label on the rightmost diagram, we conclude that 
$\cE_{6\omega_6}(-3)$ is acyclic. Proceeding in the same way with the 
other bundles, we conclude the proof. \qed

\medskip\noindent {\it Remark}.
 Observe that the irreducible vector bundle $\wedge^2\cS$ is {\it not} exceptional. 
Indeed, if $U$ is again the vector representation of $\mathrm{Spin}_{10}$, $\wedge^2U\otimes \wedge^2U$
contains $\wedge^4U$, which is an irreducible (but not fundamental) representation, contained
in the tensor product of the two half-spin representations. This implies that 
$End(\wedge^2\cS)$ contains $\cE_{\om_2+\om_3}(-2)$, which is not acyclic. Indeed, 
$s_{\alpha_1}(\om_2+\om_3-2\om_1+\rho)=\om_2+\rho$ is strictly dominant, hence
$$H^1(\OP2, \cE_{\om_2+\om_3}(-2))=\fe_6.$$

\subsection{A maximal exceptional sequence}

Recall that an exceptional sequence of sheaves on a projective variety $X$ is 
a sequence $\cF_1,\ldots , \cF_m$ of exceptional sheaves such that 
$$Ext^q(\cF_i,\cF_j)=0 \qquad \forall q\ge 0, \quad \forall i>j.$$
It is strongly exceptional if moreover
$$Ext^q(\cF_i,\cF_j)=0 \qquad \forall q>0, \quad \forall i\le j.$$
Since $\OP2$ has index $12$, it follows from the Kodaira vanishing theorem that
the sequence 
$$\cO_{\OP2}, \cO_{\OP2}(1), \ldots,   \cO_{\OP2}(10), \cO_{\OP2}(11)$$
is strongly exceptional. On the other hand, it is easy to see that the classes in K-theory
of the members of an exceptional sequence are  linearly independent (see \cite{bo}). For rational
homogeneous spaces, the K-theory is a free $\ZZ$-module admitting for basis 
the classes of the structure sheaves of the Schubert varieties. The 
length of a maximal exceptional sequence is expected to coincide with the rank 
of the K-theory, that is, the number of Schubert classes, which is also the 
topological Euler characteristic of the variety. For the Cayley plane this 
number is equal to $27$, so we expect to be able to enlarge the preceeding 
exceptional sequence of line bundles. For this we will use the exceptional 
bundles $\cS, \cS_2$ and $\cS_3$, and will apply Bott's theorem
again and again. 

\begin{lemm}\label{l1}
The bundle $\cS(-i)$ is acyclic for $1\le i\le 12$.  
\end{lemm}

\proof We play the same game as above, starting with the weight $\omega_6-i\omega_1+\rho$. 
At each step, the weight we get either has a zero coefficient, in which case the game stops 
and we conclude that we started with a singular weight, or there is a negative coefficient 
and we apply the corresponding simple reflexion. This goes as follows:

\medskip
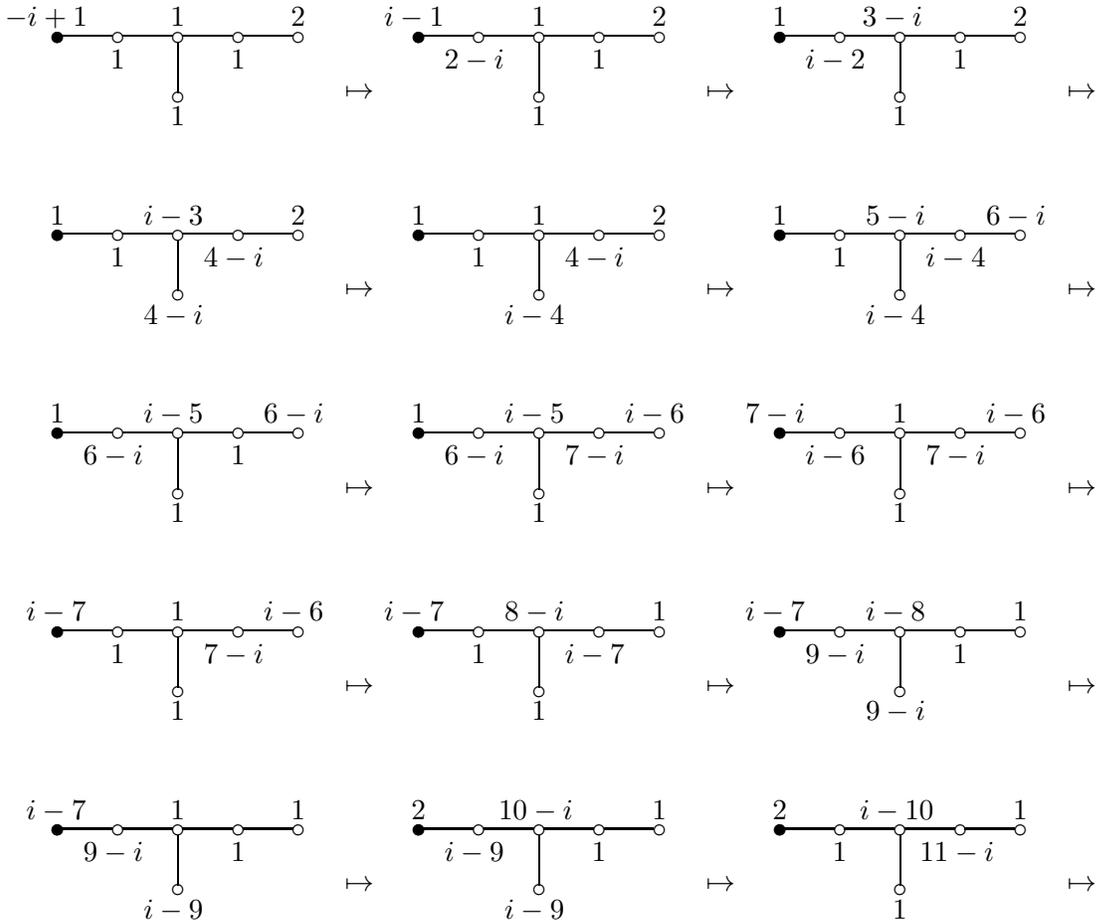
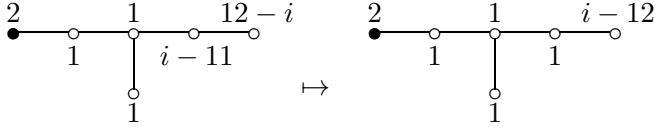
\begin{center}
\setlength{\unitlength}{4mm}
\begin{picture}(20,5)(7,0)
\multiput(-.3,3.8)(2,0){5}{$\circ$}
\multiput(0.1,4.08)(2,0){4}{\line(1,0){1.65}}
\put(3.7,1.8){$\circ$}
\put(-.3,3.8){$\bullet$}
\put(-1.8,4.4){$-i+1$}
\put(1.7,3){$1$}
\put(3.7,4.4){$1$}
\put(5.7,3){$1$}
\put(7.7,4.4){$2$}
\put(3.7,1.1){$1$}
\put(3.94,2.18){\line(0,1){1.68}}

\put(9.5,2){$\mapsto$}

\multiput(11.7,3.8)(2,0){5}{$\circ$}
\multiput(12.1,4.08)(2,0){4}{\line(1,0){1.65}}
\put(15.7,1.8){$\circ$}
\put(11.7,3.8){$\bullet$}
\put(10.8,4.4){$i-1$}
\put(12.8,3){$2-i$}
\put(15.7,4.4){$1$}
\put(17.7,3){$1$}
\put(19.7,4.4){$2$}
\put(15.7,1.1){$1$}
\put(15.94,2.18){\line(0,1){1.68}}

\put(21.5,2){$\mapsto$}

\multiput(23.7,3.8)(2,0){5}{$\circ$}
\multiput(24.1,4.08)(2,0){4}{\line(1,0){1.65}}
\put(27.7,1.8){$\circ$}
\put(23.7,3.8){$\bullet$}
\put(23.7,4.4){$1$}
\put(24.8,3){$i-2$}
\put(26.7,4.4){$3-i$}
\put(29.7,3){$1$}
\put(31.7,4.4){$2$}
\put(27.7,1.1){$1$}
\put(27.94,2.18){\line(0,1){1.68}}

\put(33.5,2){$\mapsto$}
\end{picture}
\end{center}
 
\medskip
\begin{center}
\setlength{\unitlength}{4mm}
\begin{picture}(20,5)(7,0)
\multiput(-.3,3.8)(2,0){5}{$\circ$}
\multiput(0.1,4.08)(2,0){4}{\line(1,0){1.65}}
\put(3.7,1.8){$\circ$}
\put(-.3,3.8){$\bullet$}
\put(-.3,4.4){$1$}
\put(1.7,3){$1$}
\put(2.8,4.4){$i-3$}
\put(4.8,3){$4-i$}
\put(7.7,4.4){$2$}
\put(2.8,1.1){$4-i$}
\put(3.94,2.18){\line(0,1){1.68}}

\put(9.5,2){$\mapsto$}

\multiput(11.7,3.8)(2,0){5}{$\circ$}
\multiput(12.1,4.08)(2,0){4}{\line(1,0){1.65}}
\put(15.7,1.8){$\circ$}
\put(11.7,3.8){$\bullet$}
\put(11.7,4.4){$1$}
\put(13.7,3){$1$}
\put(15.7,4.4){$1$}
\put(16.8,3){$4-i$}
\put(19.7,4.4){$2$}
\put(14.8,1.1){$i-4$}
\put(15.94,2.18){\line(0,1){1.68}}

\put(21.5,2){$\mapsto$}

\multiput(23.7,3.8)(2,0){5}{$\circ$}
\multiput(24.1,4.08)(2,0){4}{\line(1,0){1.65}}
\put(27.7,1.8){$\circ$}
\put(23.7,3.8){$\bullet$}
\put(23.7,4.4){$1$}
\put(25.7,3){$1$}
\put(26.8,4.4){$5-i$}
\put(28.8,3){$i-4$}
\put(30.8,4.4){$6-i$}
\put(26.8,1.1){$i-4$}
\put(27.94,2.18){\line(0,1){1.68}}

\put(33.5,2){$\mapsto$}
\end{picture}
\end{center}

\medskip
\begin{center}
\setlength{\unitlength}{4mm}
\begin{picture}(20,5)(7,0)
\multiput(-.3,3.8)(2,0){5}{$\circ$}
\multiput(0.1,4.08)(2,0){4}{\line(1,0){1.65}}
\put(3.7,1.8){$\circ$}
\put(-.3,3.8){$\bullet$}
\put(-.3,4.4){$1$}
\put(.8,3){$6-i$}
\put(2.8,4.4){$i-5$}
\put(5.7,3){$1$}
\put(6.8,4.4){$6-i$}
\put(3.7,1.1){$1$}
\put(3.94,2.18){\line(0,1){1.68}}

\put(9.5,2){$\mapsto$}

\multiput(11.7,3.8)(2,0){5}{$\circ$}
\multiput(12.1,4.08)(2,0){4}{\line(1,0){1.65}}
\put(15.7,1.8){$\circ$}
\put(11.7,3.8){$\bullet$}
\put(11.7,4.4){$1$}
\put(12.8,3){$6-i$}
\put(14.8,4.4){$i-5$}
\put(16.8,3){$7-i$}
\put(18.8,4.4){$i-6$}
\put(15.7,1.1){$1$}
\put(15.94,2.18){\line(0,1){1.68}}

\put(21.5,2){$\mapsto$}

\multiput(23.7,3.8)(2,0){5}{$\circ$}
\multiput(24.1,4.08)(2,0){4}{\line(1,0){1.65}}
\put(27.7,1.8){$\circ$}
\put(23.7,3.8){$\bullet$}
\put(22.8,4.4){$7-i$}
\put(24.8,3){$i-6$}
\put(27.7,4.4){$1$}
\put(28.8,3){$7-i$}
\put(30.8,4.4){$i-6$}
\put(27.7,1.1){$1$}
\put(27.94,2.18){\line(0,1){1.68}}

\put(33.5,2){$\mapsto$}
\end{picture}
\end{center}

\medskip
\begin{center}
\setlength{\unitlength}{4mm}
\begin{picture}(20,5)(7,0)
\multiput(-.3,3.8)(2,0){5}{$\circ$}
\multiput(0.1,4.08)(2,0){4}{\line(1,0){1.65}}
\put(3.7,1.8){$\circ$}
\put(-.3,3.8){$\bullet$}
\put(-1.1,4.4){$i-7$}
\put(1.7,3){$1$}
\put(3.7,4.4){$1$}
\put(4.8,3){$7-i$}
\put(6.8,4.4){$i-6$}
\put(3.7,1.1){$1$}
\put(3.94,2.18){\line(0,1){1.68}}

\put(9.5,2){$\mapsto$}

\multiput(11.7,3.8)(2,0){5}{$\circ$}
\multiput(12.1,4.08)(2,0){4}{\line(1,0){1.65}}
\put(15.7,1.8){$\circ$}
\put(11.7,3.8){$\bullet$}
\put(10.8,4.4){$i-7$}
\put(13.7,3){$1$}
\put(14.8,4.4){$8-i$}
\put(16.8,3){$i-7$}
\put(19.7,4.4){$1$}
\put(15.7,1.1){$1$}
\put(15.94,2.18){\line(0,1){1.68}}

\put(21.5,2){$\mapsto$}

\multiput(23.7,3.8)(2,0){5}{$\circ$}
\multiput(24.1,4.08)(2,0){4}{\line(1,0){1.65}}
\put(27.7,1.8){$\circ$}
\put(23.7,3.8){$\bullet$}
\put(22.8,4.4){$i-7$}
\put(24.8,3){$9-i$}
\put(26.8,4.4){$i-8$}
\put(29.7,3){$1$}
\put(31.7,4.4){$1$}
\put(26.8,1.1){$9-i$}
\put(27.94,2.18){\line(0,1){1.68}}

\put(33.5,2){$\mapsto$}
\end{picture}
\end{center}

\medskip
\begin{center}
\setlength{\unitlength}{4mm}
\begin{picture}(20,5)(7,0)
\multiput(-.3,3.8)(2,0){5}{$\circ$}
\multiput(0.1,4.08)(2,0){4}{\line(1,0){1.65}}
\put(3.7,1.8){$\circ$}
\put(-.3,3.8){$\bullet$}
\put(-1.1,4.4){$i-7$}
\put(.8,3){$9-i$}
\put(3.7,4.4){$1$}
\put(5.7,3){$1$}
\put(7.7,4.4){$1$}
\put(2.8,1.1){$i-9$}
\put(3.94,2.18){\line(0,1){1.68}}

\put(9.5,2){$\mapsto$}

\multiput(11.7,3.8)(2,0){5}{$\circ$}
\multiput(12.1,4.08)(2,0){4}{\line(1,0){1.65}}
\put(15.7,1.8){$\circ$}
\put(11.7,3.8){$\bullet$}
\put(11.7,4.4){$2$}
\put(12.8,3){$i-9$}
\put(14.6,4.4){$10-i$}
\put(17.7,3){$1$}
\put(19.7,4.4){$1$}
\put(14.8,1.1){$i-9$}
\put(15.94,2.18){\line(0,1){1.68}}

\put(21.5,2){$\mapsto$}

\multiput(23.7,3.8)(2,0){5}{$\circ$}
\multiput(24.1,4.08)(2,0){4}{\line(1,0){1.65}}
\put(27.7,1.8){$\circ$}
\put(23.7,3.8){$\bullet$}
\put(23.7,4.4){$2$}
\put(25.7,3){$1$}
\put(26.6,4.4){$i-10$}
\put(28.6,3){$11-i$}
\put(31.7,4.4){$1$}
\put(27.7,1.1){$1$}
\put(27.94,2.18){\line(0,1){1.68}}

\put(33.5,2){$\mapsto$}
\end{picture}
\end{center}

\medskip
\begin{center}
\setlength{\unitlength}{4mm}
\begin{picture}(20,5)(7,0)
\multiput(-.3,3.8)(2,0){5}{$\circ$}
\multiput(0.1,4.08)(2,0){4}{\line(1,0){1.65}}
\put(3.7,1.8){$\circ$}
\put(-.3,3.8){$\bullet$}
\put(-.3,4.4){$2$}
\put(1.7,3){$1$}
\put(3.7,4.4){$1$}
\put(4.8,3){$i-11$}
\put(6.8,4.4){$12-i$}
\put(3.7,1.1){$1$}
\put(3.94,2.18){\line(0,1){1.68}}

\put(9.5,2){$\mapsto$}

\multiput(11.7,3.8)(2,0){5}{$\circ$}
\multiput(12.1,4.08)(2,0){4}{\line(1,0){1.65}}
\put(15.7,1.8){$\circ$}
\put(11.7,3.8){$\bullet$}
\put(11.7,4.4){$2$}
\put(13.7,3){$1$}
\put(15.7,4.4){$1$}
\put(17.7,3){$1$}
\put(18.8,4.4){$i-12$}
\put(15.7,1.1){$1$}
\put(15.94,2.18){\line(0,1){1.68}}
\end{picture}
\end{center}
This concludes the proof. \qed

\medskip
Note that for $i=13$ we finally get the strictly dominant weight $\om_1+\rho$. 
Since we needed to apply $16$ simple reflexions, we conclude by Bott's theorem that 
$$H^{16}(\OP2,\cS(-13))=V_{\om_1}^\vee.$$
But by Serre duality,  $H^{16}(\OP2,\cS(-13))$ is dual to $H^0(\OP2,\cS^\vee(1))=H^0(\OP2,\cS)$
which, by Borel-Weil, is $V_{\om_6}^\vee\simeq V_{\om_1}$. This is a way to check that the 
computation above, and
those of the same type that will follow, are indeed correct. 

\smallskip
The same statement holds for our two other exceptional bundles:

\begin{lemm}\label{l2}
$\cS_2(-i)$ and $\cS_3(-i)$ are acyclic for $1\le i\le 12$.
\end{lemm}

Now consider their endomorphism bundles:

\begin{lemm}\label{l3}
$End(\cS)(-i)$ is acyclic for $1\le i\le 11$.  
\end{lemm}

\proof We have seen that $End(\cS)=\cS_2(-1)\oplus\cO_{\OP2}\oplus \cE_{\omega_5}(-1)$. 
We already know that $\cS_2(-i-1)$ and $\cO_{\OP2}(-i)$ are acyclic for $1\le i\le 11$.
There remains to treat the case of $\cE_{\omega_5}(-i-1)$, which we do as above.
After adding $\rho$ to $\omega_5-(i+1)\om_1$, we get successively:

\medskip
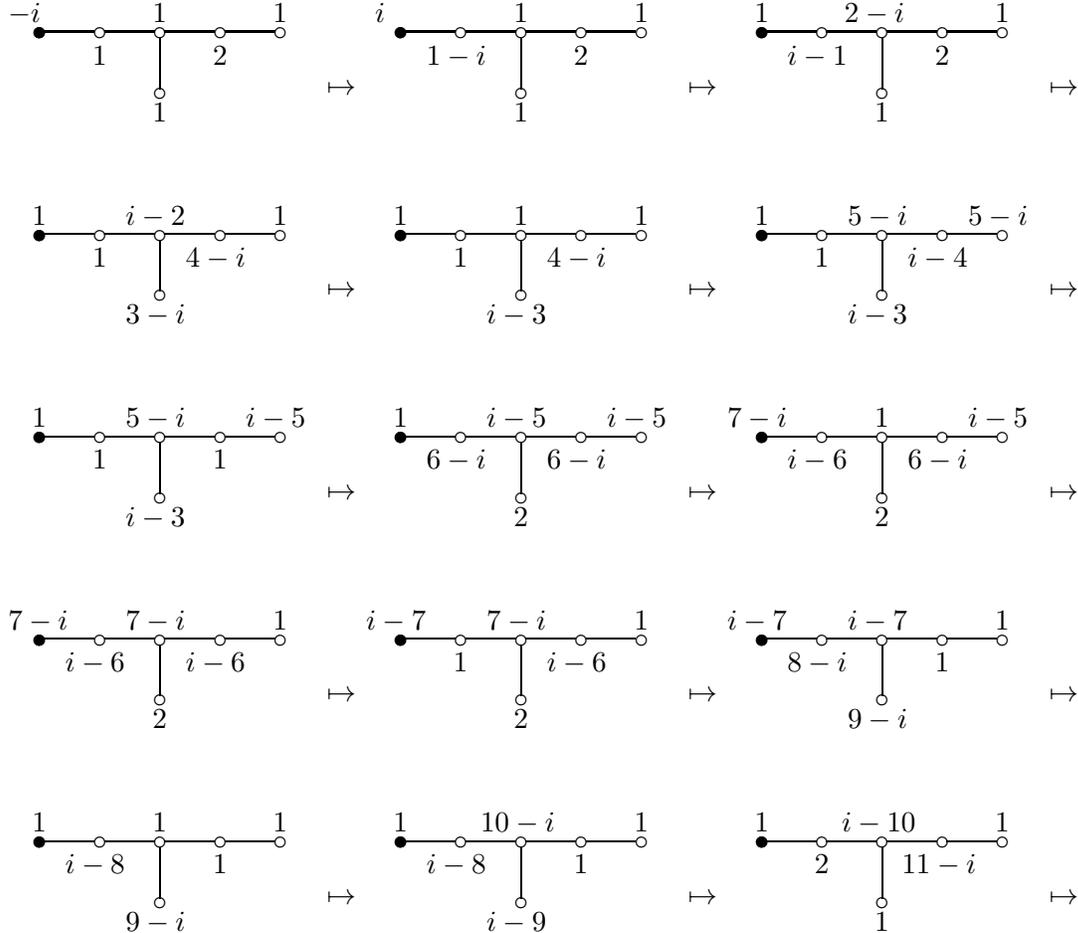
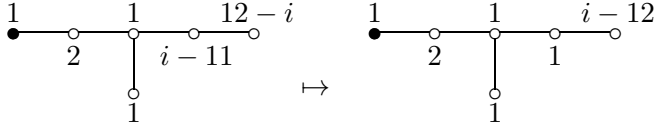
\begin{center}
\setlength{\unitlength}{4mm}
\begin{picture}(20,5)(7,0)
\multiput(-.3,3.8)(2,0){5}{$\circ$}
\multiput(0.1,4.08)(2,0){4}{\line(1,0){1.65}}
\put(3.7,1.8){$\circ$}
\put(-.3,3.8){$\bullet$}
\put(-1.1,4.4){$-i$}
\put(1.7,3){$1$}
\put(3.7,4.4){$1$}
\put(5.7,3){$2$}
\put(7.7,4.4){$1$}
\put(3.7,1.1){$1$}
\put(3.94,2.18){\line(0,1){1.68}}

\put(9.5,2){$\mapsto$}

\multiput(11.7,3.8)(2,0){5}{$\circ$}
\multiput(12.1,4.08)(2,0){4}{\line(1,0){1.65}}
\put(15.7,1.8){$\circ$}
\put(11.7,3.8){$\bullet$}
\put(11.1,4.4){$i$}
\put(12.8,3){$1-i$}
\put(15.7,4.4){$1$}
\put(17.7,3){$2$}
\put(19.7,4.4){$1$}
\put(15.7,1.1){$1$}
\put(15.94,2.18){\line(0,1){1.68}}

\put(21.5,2){$\mapsto$}

\multiput(23.7,3.8)(2,0){5}{$\circ$}
\multiput(24.1,4.08)(2,0){4}{\line(1,0){1.65}}
\put(27.7,1.8){$\circ$}
\put(23.7,3.8){$\bullet$}
\put(23.7,4.4){$1$}
\put(24.8,3){$i-1$}
\put(26.7,4.4){$2-i$}
\put(29.7,3){$2$}
\put(31.7,4.4){$1$}
\put(27.7,1.1){$1$}
\put(27.94,2.18){\line(0,1){1.68}}

\put(33.5,2){$\mapsto$}
\end{picture}
\end{center}
 
\medskip
\begin{center}
\setlength{\unitlength}{4mm}
\begin{picture}(20,5)(7,0)
\multiput(-.3,3.8)(2,0){5}{$\circ$}
\multiput(0.1,4.08)(2,0){4}{\line(1,0){1.65}}
\put(3.7,1.8){$\circ$}
\put(-.3,3.8){$\bullet$}
\put(-.3,4.4){$1$}
\put(1.7,3){$1$}
\put(2.8,4.4){$i-2$}
\put(4.8,3){$4-i$}
\put(7.7,4.4){$1$}
\put(2.8,1.1){$3-i$}
\put(3.94,2.18){\line(0,1){1.68}}

\put(9.5,2){$\mapsto$}

\multiput(11.7,3.8)(2,0){5}{$\circ$}
\multiput(12.1,4.08)(2,0){4}{\line(1,0){1.65}}
\put(15.7,1.8){$\circ$}
\put(11.7,3.8){$\bullet$}
\put(11.7,4.4){$1$}
\put(13.7,3){$1$}
\put(15.7,4.4){$1$}
\put(16.8,3){$4-i$}
\put(19.7,4.4){$1$}
\put(14.8,1.1){$i-3$}
\put(15.94,2.18){\line(0,1){1.68}}

\put(21.5,2){$\mapsto$}

\multiput(23.7,3.8)(2,0){5}{$\circ$}
\multiput(24.1,4.08)(2,0){4}{\line(1,0){1.65}}
\put(27.7,1.8){$\circ$}
\put(23.7,3.8){$\bullet$}
\put(23.7,4.4){$1$}
\put(25.7,3){$1$}
\put(26.8,4.4){$5-i$}
\put(28.8,3){$i-4$}
\put(30.8,4.4){$5-i$}
\put(26.8,1.1){$i-3$}
\put(27.94,2.18){\line(0,1){1.68}}

\put(33.5,2){$\mapsto$}
\end{picture}
\end{center}

\medskip
\begin{center}
\setlength{\unitlength}{4mm}
\begin{picture}(20,5)(7,0)
\multiput(-.3,3.8)(2,0){5}{$\circ$}
\multiput(0.1,4.08)(2,0){4}{\line(1,0){1.65}}
\put(3.7,1.8){$\circ$}
\put(-.3,3.8){$\bullet$}
\put(-.3,4.4){$1$}
\put(1.7,3){$1$}
\put(2.8,4.4){$5-i$}
\put(5.7,3){$1$}
\put(6.8,4.4){$i-5$}
\put(2.8,1.1){$i-3$}
\put(3.94,2.18){\line(0,1){1.68}}

\put(9.5,2){$\mapsto$}

\multiput(11.7,3.8)(2,0){5}{$\circ$}
\multiput(12.1,4.08)(2,0){4}{\line(1,0){1.65}}
\put(15.7,1.8){$\circ$}
\put(11.7,3.8){$\bullet$}
\put(11.7,4.4){$1$}
\put(12.8,3){$6-i$}
\put(14.8,4.4){$i-5$}
\put(16.8,3){$6-i$}
\put(18.8,4.4){$i-5$}
\put(15.7,1.1){$2$}
\put(15.94,2.18){\line(0,1){1.68}}

\put(21.5,2){$\mapsto$}

\multiput(23.7,3.8)(2,0){5}{$\circ$}
\multiput(24.1,4.08)(2,0){4}{\line(1,0){1.65}}
\put(27.7,1.8){$\circ$}
\put(23.7,3.8){$\bullet$}
\put(22.8,4.4){$7-i$}
\put(24.8,3){$i-6$}
\put(27.7,4.4){$1$}
\put(28.8,3){$6-i$}
\put(30.8,4.4){$i-5$}
\put(27.7,1.1){$2$}
\put(27.94,2.18){\line(0,1){1.68}}

\put(33.5,2){$\mapsto$}
\end{picture}
\end{center}

\medskip
\begin{center}
\setlength{\unitlength}{4mm}
\begin{picture}(20,5)(7,0)
\multiput(-.3,3.8)(2,0){5}{$\circ$}
\multiput(0.1,4.08)(2,0){4}{\line(1,0){1.65}}
\put(3.7,1.8){$\circ$}
\put(-.3,3.8){$\bullet$}
\put(-1.1,4.4){$7-i$}
\put(.8,3){$i-6$}
\put(2.8,4.4){$7-i$}
\put(4.8,3){$i-6$}
\put(7.7,4.4){$1$}
\put(3.7,1.1){$2$}
\put(3.94,2.18){\line(0,1){1.68}}

\put(9.5,2){$\mapsto$}

\multiput(11.7,3.8)(2,0){5}{$\circ$}
\multiput(12.1,4.08)(2,0){4}{\line(1,0){1.65}}
\put(15.7,1.8){$\circ$}
\put(11.7,3.8){$\bullet$}
\put(10.8,4.4){$i-7$}
\put(13.7,3){$1$}
\put(14.8,4.4){$7-i$}
\put(16.8,3){$i-6$}
\put(19.7,4.4){$1$}
\put(15.7,1.1){$2$}
\put(15.94,2.18){\line(0,1){1.68}}

\put(21.5,2){$\mapsto$}

\multiput(23.7,3.8)(2,0){5}{$\circ$}
\multiput(24.1,4.08)(2,0){4}{\line(1,0){1.65}}
\put(27.7,1.8){$\circ$}
\put(23.7,3.8){$\bullet$}
\put(22.8,4.4){$i-7$}
\put(24.8,3){$8-i$}
\put(26.8,4.4){$i-7$}
\put(29.7,3){$1$}
\put(31.7,4.4){$1$}
\put(26.8,1.1){$9-i$}
\put(27.94,2.18){\line(0,1){1.68}}

\put(33.5,2){$\mapsto$}
\end{picture}
\end{center}

\medskip
\begin{center}
\setlength{\unitlength}{4mm}
\begin{picture}(20,5)(7,0)
\multiput(-.3,3.8)(2,0){5}{$\circ$}
\multiput(0.1,4.08)(2,0){4}{\line(1,0){1.65}}
\put(3.7,1.8){$\circ$}
\put(-.3,3.8){$\bullet$}
\put(-.3,4.4){$1$}
\put(.8,3){$i-8$}
\put(3.7,4.4){$1$}
\put(5.7,3){$1$}
\put(7.7,4.4){$1$}
\put(2.8,1.1){$9-i$}
\put(3.94,2.18){\line(0,1){1.68}}

\put(9.5,2){$\mapsto$}

\multiput(11.7,3.8)(2,0){5}{$\circ$}
\multiput(12.1,4.08)(2,0){4}{\line(1,0){1.65}}
\put(15.7,1.8){$\circ$}
\put(11.7,3.8){$\bullet$}
\put(11.7,4.4){$1$}
\put(12.8,3){$i-8$}
\put(14.6,4.4){$10-i$}
\put(17.7,3){$1$}
\put(19.7,4.4){$1$}
\put(14.8,1.1){$i-9$}
\put(15.94,2.18){\line(0,1){1.68}}

\put(21.5,2){$\mapsto$}

\multiput(23.7,3.8)(2,0){5}{$\circ$}
\multiput(24.1,4.08)(2,0){4}{\line(1,0){1.65}}
\put(27.7,1.8){$\circ$}
\put(23.7,3.8){$\bullet$}
\put(23.7,4.4){$1$}
\put(25.7,3){$2$}
\put(26.6,4.4){$i-10$}
\put(28.6,3){$11-i$}
\put(31.7,4.4){$1$}
\put(27.7,1.1){$1$}
\put(27.94,2.18){\line(0,1){1.68}}

\put(33.5,2){$\mapsto$}
\end{picture}
\end{center}

\medskip
\begin{center}
\setlength{\unitlength}{4mm}
\begin{picture}(20,5)(7,0)
\multiput(-.3,3.8)(2,0){5}{$\circ$}
\multiput(0.1,4.08)(2,0){4}{\line(1,0){1.65}}
\put(3.7,1.8){$\circ$}
\put(-.3,3.8){$\bullet$}
\put(-.3,4.4){$1$}
\put(1.7,3){$2$}
\put(3.7,4.4){$1$}
\put(4.8,3){$i-11$}
\put(6.8,4.4){$12-i$}
\put(3.7,1.1){$1$}
\put(3.94,2.18){\line(0,1){1.68}}

\put(9.5,2){$\mapsto$}

\multiput(11.7,3.8)(2,0){5}{$\circ$}
\multiput(12.1,4.08)(2,0){4}{\line(1,0){1.65}}
\put(15.7,1.8){$\circ$}
\put(11.7,3.8){$\bullet$}
\put(11.7,4.4){$1$}
\put(13.7,3){$2$}
\put(15.7,4.4){$1$}
\put(17.7,3){$1$}
\put(18.8,4.4){$i-12$}
\put(15.7,1.1){$1$}
\put(15.94,2.18){\line(0,1){1.68}}
\end{picture}
\end{center}
This concludes the proof. \qed

\begin{lemm}\label{l4}
$End(\cS_2)(-i)$ is acyclic for $1\le i\le 2$.  
\end{lemm}

\proof We have seen how to decompose $End(\cS_2)$ into irreducible bundles. We need 
to apply Bott's theorem to each component. Consider for example the component 
$\cE_{2\om_5}(-2)$. After twisting by $\cO_{\OP2}(-i)$ and adding $\rho$ to the corresponding weight, we get:

\medskip
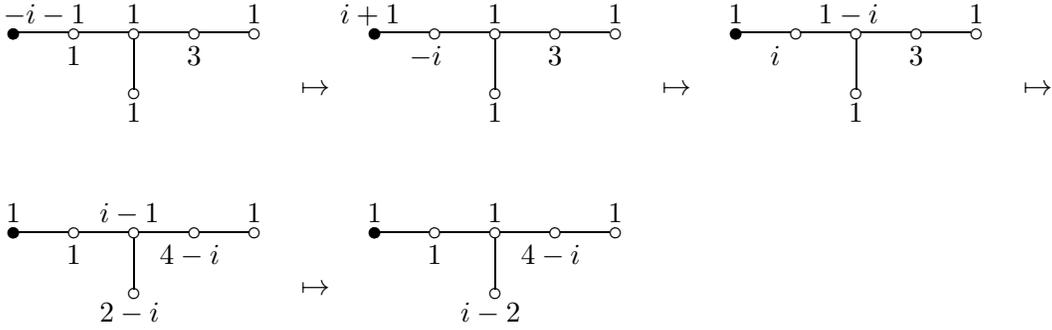
\begin{center}
\setlength{\unitlength}{4mm}
\begin{picture}(20,5)(7,0)
\multiput(-.3,3.8)(2,0){5}{$\circ$}
\multiput(0.1,4.08)(2,0){4}{\line(1,0){1.65}}
\put(3.7,1.8){$\circ$}
\put(-.3,3.8){$\bullet$}
\put(-.4,4.4){$-i-1$}
\put(1.7,3){$1$}
\put(3.7,4.4){$1$}
\put(5.7,3){$3$}
\put(7.7,4.4){$1$}
\put(3.7,1.1){$1$}
\put(3.94,2.18){\line(0,1){1.68}}

\put(9.5,2){$\mapsto$}

\multiput(11.7,3.8)(2,0){5}{$\circ$}
\multiput(12.1,4.08)(2,0){4}{\line(1,0){1.65}}
\put(15.7,1.8){$\circ$}
\put(11.7,3.8){$\bullet$}
\put(10.8,4.4){$i+1$}
\put(13.1,3){$-i$}
\put(15.7,4.4){$1$}
\put(17.7,3){$3$}
\put(19.7,4.4){$1$}
\put(15.7,1.1){$1$}
\put(15.94,2.18){\line(0,1){1.68}}

\put(21.5,2){$\mapsto$}

\multiput(23.7,3.8)(2,0){5}{$\circ$}
\multiput(24.1,4.08)(2,0){4}{\line(1,0){1.65}}
\put(27.7,1.8){$\circ$}
\put(23.7,3.8){$\bullet$}
\put(23.7,4.4){$1$}
\put(25.1,3){$i$}
\put(26.7,4.4){$1-i$}
\put(29.7,3){$3$}
\put(31.7,4.4){$1$}
\put(27.7,1.1){$1$}
\put(27.94,2.18){\line(0,1){1.68}}

\put(33.5,2){$\mapsto$}
\end{picture}
\end{center}
 
\medskip
\begin{center}
\setlength{\unitlength}{4mm}
\begin{picture}(20,5)(7,0)
\multiput(-.3,3.8)(2,0){5}{$\circ$}
\multiput(0.1,4.08)(2,0){4}{\line(1,0){1.65}}
\put(3.7,1.8){$\circ$}
\put(-.3,3.8){$\bullet$}
\put(-.3,4.4){$1$}
\put(1.7,3){$1$}
\put(2.8,4.4){$i-1$}
\put(4.8,3){$4-i$}
\put(7.7,4.4){$1$}
\put(2.8,1.1){$2-i$}
\put(3.94,2.18){\line(0,1){1.68}}

\put(9.5,2){$\mapsto$}

\multiput(11.7,3.8)(2,0){5}{$\circ$}
\multiput(12.1,4.08)(2,0){4}{\line(1,0){1.65}}
\put(15.7,1.8){$\circ$}
\put(11.7,3.8){$\bullet$}
\put(11.7,4.4){$1$}
\put(13.7,3){$1$}
\put(15.7,4.4){$1$}
\put(16.8,3){$4-i$}
\put(19.7,4.4){$1$}
\put(14.8,1.1){$i-2$}
\put(15.94,2.18){\line(0,1){1.68}}
\end{picture}
\end{center}

For $i=1,2$ we get singular weights, as claimed. But note that for $i=3$, 
the last weight above is $\rho$, so that 
that $H^3(\OP2,\cE_{2\om_5}(-5))=\CC$. In particular   $End(\cS_2)(-3)$ is not acyclic. 
Examining the other components we can easily complete the proof that $End(\cS_2)(-1)$ 
and $End(\cS_2)(-2)$ are both  acyclic.\qed

\medskip In a completely similar way, we check that:

\begin{lemm}\label{l5}
$End(\cS_3)(-1)$ is acyclic.
\end{lemm}

\begin{lemm}\label{l6}
$\cS_2\otimes\cS(-i-1)$ is acyclic for $1\le i\le 12$.  
\end{lemm}

\proof Use the decomposition, that we obtain e.g. using LiE, 
$$\cS_2\otimes\cS = \cS_3\oplus \cS(1)\oplus\cE_{\omega_5+\omega_6}.$$
The first two factors have already been considered. The third one is treated in the same way. \qed

\begin{lemm}\label{l7}
$\cS_3\otimes\cS(-i-1)$ is acyclic for $1\le i\le 6$.  
\end{lemm}

\proof Here the relevant decomposition is
$$\cS_3\otimes\cS = \cE_{4\omega_6}\oplus\cS_2(1)\oplus \cE_{\omega_5+2\omega_6}.$$
The most limiting term is the first one, since it gives rise to the sequence:

\medskip
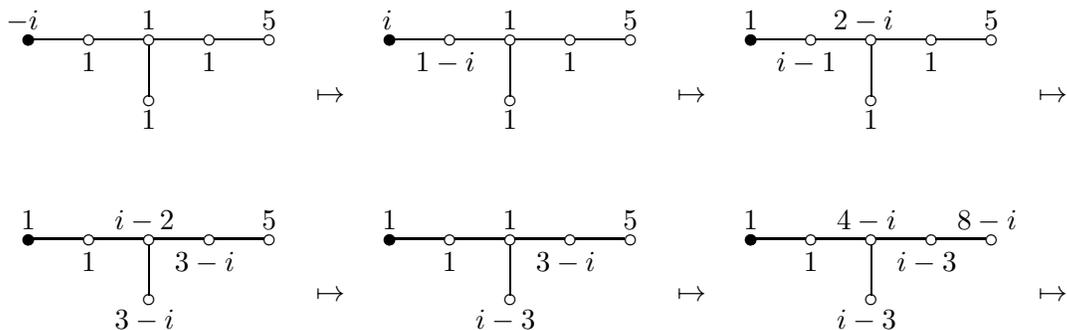
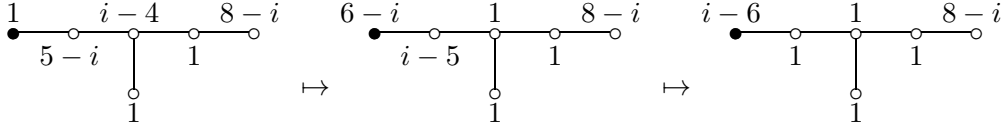
\begin{center}
\setlength{\unitlength}{4mm}
\begin{picture}(20,5)(7,0)
\multiput(-.3,3.8)(2,0){5}{$\circ$}
\multiput(0.1,4.08)(2,0){4}{\line(1,0){1.65}}
\put(3.7,1.8){$\circ$}
\put(-.3,3.8){$\bullet$}
\put(-.8,4.4){$-i$}
\put(1.7,3){$1$}
\put(3.7,4.4){$1$}
\put(5.7,3){$1$}
\put(7.7,4.4){$5$}
\put(3.7,1.1){$1$}
\put(3.94,2.18){\line(0,1){1.68}}

\put(9.5,2){$\mapsto$}

\multiput(11.7,3.8)(2,0){5}{$\circ$}
\multiput(12.1,4.08)(2,0){4}{\line(1,0){1.65}}
\put(15.7,1.8){$\circ$}
\put(11.7,3.8){$\bullet$}
\put(11.7,4.4){$i$}
\put(12.8,3){$1-i$}
\put(15.7,4.4){$1$}
\put(17.7,3){$1$}
\put(19.7,4.4){$5$}
\put(15.7,1.1){$1$}
\put(15.94,2.18){\line(0,1){1.68}}

\put(21.5,2){$\mapsto$}

\multiput(23.7,3.8)(2,0){5}{$\circ$}
\multiput(24.1,4.08)(2,0){4}{\line(1,0){1.65}}
\put(27.7,1.8){$\circ$}
\put(23.7,3.8){$\bullet$}
\put(23.7,4.4){$1$}
\put(24.8,3){$i-1$}
\put(26.7,4.4){$2-i$}
\put(29.7,3){$1$}
\put(31.7,4.4){$5$}
\put(27.7,1.1){$1$}
\put(27.94,2.18){\line(0,1){1.68}}

\put(33.5,2){$\mapsto$}
\end{picture}
\end{center}
 
\medskip
\begin{center}
\setlength{\unitlength}{4mm}
\begin{picture}(20,5)(7,0)
\multiput(-.3,3.8)(2,0){5}{$\circ$}
\multiput(0.1,4.08)(2,0){4}{\line(1,0){1.65}}
\put(3.7,1.8){$\circ$}
\put(-.3,3.8){$\bullet$}
\put(-.3,4.4){$1$}
\put(1.7,3){$1$}
\put(2.8,4.4){$i-2$}
\put(4.8,3){$3-i$}
\put(7.7,4.4){$5$}
\put(2.8,1.1){$3-i$}
\put(3.94,2.18){\line(0,1){1.68}}

\put(9.5,2){$\mapsto$}

\multiput(11.7,3.8)(2,0){5}{$\circ$}
\multiput(12.1,4.08)(2,0){4}{\line(1,0){1.65}}
\put(15.7,1.8){$\circ$}
\put(11.7,3.8){$\bullet$}
\put(11.7,4.4){$1$}
\put(13.7,3){$1$}
\put(15.7,4.4){$1$}
\put(16.8,3){$3-i$}
\put(19.7,4.4){$5$}
\put(14.8,1.1){$i-3$}
\put(15.94,2.18){\line(0,1){1.68}}

\put(21.5,2){$\mapsto$}

\multiput(23.7,3.8)(2,0){5}{$\circ$}
\multiput(24.1,4.08)(2,0){4}{\line(1,0){1.65}}
\put(27.7,1.8){$\circ$}
\put(23.7,3.8){$\bullet$}
\put(23.7,4.4){$1$}
\put(25.7,3){$1$}
\put(26.8,4.4){$4-i$}
\put(28.8,3){$i-3$}
\put(30.8,4.4){$8-i$}
\put(26.8,1.1){$i-3$}
\put(27.94,2.18){\line(0,1){1.68}}

\put(33.5,2){$\mapsto$}
\end{picture}
\end{center}

\medskip
\begin{center}
\setlength{\unitlength}{4mm}
\begin{picture}(20,5)(7,0)
\multiput(-.3,3.8)(2,0){5}{$\circ$}
\multiput(0.1,4.08)(2,0){4}{\line(1,0){1.65}}
\put(3.7,1.8){$\circ$}
\put(-.3,3.8){$\bullet$}
\put(-.3,4.4){$1$}
\put(.8,3){$5-i$}
\put(2.8,4.4){$i-4$}
\put(5.7,3){$1$}
\put(6.8,4.4){$8-i$}
\put(3.7,1.1){$1$}
\put(3.94,2.18){\line(0,1){1.68}}

\put(9.5,2){$\mapsto$}

\multiput(11.7,3.8)(2,0){5}{$\circ$}
\multiput(12.1,4.08)(2,0){4}{\line(1,0){1.65}}
\put(15.7,1.8){$\circ$}
\put(11.7,3.8){$\bullet$}
\put(10.8,4.4){$6-i$}
\put(12.8,3){$i-5$}
\put(15.7,4.4){$1$}
\put(17.7,3){$1$}
\put(18.8,4.4){$8-i$}
\put(15.7,1.1){$1$}
\put(15.94,2.18){\line(0,1){1.68}}

\put(21.5,2){$\mapsto$}

\multiput(23.7,3.8)(2,0){5}{$\circ$}
\multiput(24.1,4.08)(2,0){4}{\line(1,0){1.65}}
\put(27.7,1.8){$\circ$}
\put(23.7,3.8){$\bullet$}
\put(22.8,4.4){$i-6$}
\put(25.7,3){$1$}
\put(27.7,4.4){$1$}
\put(29.7,3){$1$}
\put(30.8,4.4){$8-i$}
\put(27.7,1.1){$1$}
\put(27.94,2.18){\line(0,1){1.68}}

\end{picture}
\end{center}
 
For $i=1,\ldots, 6$ we get singular weights, as claimed, but for $i=7$ the last weight
above is $\rho$. We can therefore conclude that $H^8(\OP2,\cE_{\omega_5+2\omega_6}(-8))=\CC$.
Therefore $\cS_3\otimes\cS(-8)$ is not acyclic. 
We conclude the proof by checking the last component. \qed

\begin{lemm}\label{l8}
$\cS_3\otimes\cS_2(-i-2)$ is acyclic for $1\le i\le 2$.  
\end{lemm}

\proof We have the decomposition:
$$\cS_3\otimes\cS_2 = \cE_{5\omega_6}\oplus \cE_{\omega_5+3\omega_6}
\oplus \cE_{2\omega_5+\omega_6}\oplus\cE_{\omega_5+\omega_6}(1)
\oplus\cS_3(1)\oplus \cS(2).$$
The last three terms have already been considered. Among the first three, the most 
limiting one is the third one, which contributes non trivially for $i=3$. But for 
$i=1,2$ all the factors are acyclic. \qed 

\medskip We can now prove our main result. 

\begin{theo}
The following sequence, of length $27$, of vector bundles on the Cayley plane $\OO\PP^2$,
\begin{eqnarray*}
 &\cO, \cS, \cO(1), \cS(1), \cO(2), \cS(2), \cO(3), \cS(3), \cO(4), \\
 &\cS_2(3), \cS(4), \cS_3(3), \cO(5), 
\cS_2(4), \cS(5), \cS_3(4), \cO(6), \cS_2(5), \\
 &\cS(6), \cO(7), \cS(7), \cO(8), \cS(8), \cO(9), \cS(9), \cO(10), \cO(11)
\end{eqnarray*}
is a maximal strongly exceptional collection.
\end{theo}

\proof This follows from the previous lemmas. Start with the exceptional 
collection $\cO,\ldots , \cO(11)$. By Lemma \ref{l3}, $\cS,\cS(1),\ldots ,\cS(9)$ 
is also an exceptional collection. According to Lemma \ref{l1}, we have 
$Hom(\cO(i),\cS(j))=0$ for $j<i\le j+12$. Moreover, since $\cS^\vee=\cS(-1)$, 
$Hom(\cS(j),\cO(i))=Hom(\cO(j+1),\cS(i))=0$ for $i\le j\le i+11$. 
This implies that the sequence  
\begin{eqnarray*}
 &\cO, \cS, \cO(1), \cS(1), \cO(2), \cS(2), \cO(3), \cS(3), \cO(4), \cS(4), \cO(5),  \\
 &\cS(5), \cO(6),\cS(6), \cO(7), \cS(7), \cO(8), \cS(8), \cO(9), \cS(9), \cO(10), \cO(11)
\end{eqnarray*}
is an exceptional collection. On the other hand, Lemmas \ref{l1}, \ref{l4}, \ref{l5} and \ref{l8} imply that 
$$\cS_2(3), \cS_3(3), \cS_2(4), \cS_3(4), \cS_2(5)$$
is also an exceptional collection. There remains to ``insert'' this collection inside the 
previous one. The compatibility conditions are the following. For $\cS_2$, Lemmas \ref{l2}, \ref{l6} 
and the fact that $\cS_2^\vee=\cS_2(-2)$ imply that we must respect the orderings 
$\cO(k)\cdots \cS_2(j)\cdots\cO(i)$ and 
$\cS(k)\cdots \cS_2(j)\cdots\cS(i)$ with $j-10\le k\le j+1$ and $j+1\le i\le j+12$.  
Concerning $\cS_3$, Lemmas \ref{l2}, \ref{l7} 
and the fact that $\cS_3^\vee=\cS_3(-3)$ imply that we must respect the orderings 
$\cO(k)\cdots \cS_3(j)\cdots\cO(i)$  with $j-9\le k\le j+2$ and $j+1\le i\le j+12$, 
and $\cS(k)\cdots \cS_3(j)\cdots\cS(i)$ with $j-4\le k\le j+1$ and $j+1\le i\le j+6$. 
The collection of the theorem is compatible with these requirements. \qed

Finally the fact that this collection is strongly exceptional is quite straightforward. Indeed,
if $i<j$ and $E_i, E_j$ are the corresponding bundles of the collection, then in most cases 
$End(E_i,E_j)$ decomposes as a sum of irreducible vector bundles $\cE_\om$ defined by a highest 
weight $\om$ which is dominant. In this case it is an immediate consequence of Bott's theorem 
that the higher 
cohomology groups vanish. Another possibility is that $\om$ has coefficient $-1$ on $\omega_1$, 
and then $\cE_\om$ is acyclic. The remaining cases are only of three types, $E_i=\cS(k)$ and  
$E_j=\cS_3(k-1)$, or $E_i=\cS_2(k)$ and  $E_j=\cS_3(k)$, or $E_i=\cS_3(k)$ and  $E_j=\cS_2(k+1)$. 
In these cases the coefficient of $\om$ on $\om_1$ can be $-2$, but then the coefficient on $\om_3$ 
is zero, and the acyclicity follows immediately. \qed

\medskip
Of course we expect this maximal exceptional collection to be full, i.e. to generate
the derived category of the Cayley plane. This would follow from Conjecture 9.1 and its Corollary 9.3
 in \cite{ku3}, but we have not been able to prove it.

A possible strategy would be to find a covering family of subvarieties, possibly of
small codimension,  for which we already have a good understanding of the derived 
category. This was the strategy used in \cite{ku1} for an inductive treatment of Grassmannians
of lines. In the Cayley plane, there are at least two natural candidates. The first one is
the family of $\OO$-lines, that is, of eight-dimensional quadrics parametrized by the dual
Cayley plane. The other one is the family of copies of the spinor variety of $\mathrm{Spin}_{10}$ 
in $\OP2$. Indeed, the union of lines in the Cayley plane passing through a given point 
is known to be a cone over this spinor variety \cite{LM}, which we can recover by taking 
hyperplane sections not containing the given point. These spinor varieties have codimension
six, and their derived category is described in \cite[6.2]{ku2}. But in both cases the codimension 
is already sufficiently big to make this strategy difficult to implement concretely.

\medskip\noindent 
{\small Laurent MANIVEL, Institut Fourier, Laboratoire de Math\'ematiques, 
UMR 5582 UJF-CNRS, BP 74, 38402 St Martin d'H\`eres Cedex, France.

\noindent E-mail : {\tt Laurent.Manivel@ujf-grenoble.fr}}

\end{document}